\documentclass[12pt]{amsart}
\usepackage{amssymb}
\usepackage{amsfonts}
\usepackage{latexsym}
\usepackage{amscd}

\newcounter{cs}
\stepcounter{cs}
\newcounter{ds}
\stepcounter{ds}

\newcommand{\casos}{\begin{itemize}}
\newcommand{\fcasos}{\end{itemize}\setcounter{cs}{1}}

\newcommand{\subcasos}{\begin{itemize}}
\newcommand{\fsubcasos}{\end{itemize}\setcounter{ds}{1}}

\newcommand{\matriu}[1]{\left(\begin{array}{#1}}
\newcommand{\fmatriu}{\end{array}\right)}

\newcommand{\ld}{{\rm lim}_{_{\kern-14pt\longrightarrow \kern3pt}}}

\newcommand{\C}{$C^*-$algebra}
\newcommand{\Cs}{$C^*-$algebras}

\newcommand{\cl}{{\rm cl}}
\newcommand{\ccr}{{\rm cr}}

\newcommand{\fl}{\rightarrow}

\newcommand{\ol}{\overline}

\newcommand{\mulr}{\mathcal{M}(R)}

\newtheorem{lem}{Lemma}[section]
\newtheorem{corol}[lem]{Corollary}
\newtheorem{theor}[lem]{Theorem}
\newtheorem{prop}[lem]{Proposition}
\newtheorem{rema}[lem]{Remark}

\newtheorem{defis}[lem]{Definitions}

\begin{document}
\title[Extensions and PullBacks in $QB-$rings]{Extensions and PullBacks in \boldmath
$QB-$rings}
\author{Pere Ara, Gert K. Pedersen and Francesc Perera}
\thanks{Partially supported by the DGESIC (Spain), the Comissionat per
Universitats i Recerca de la Generalitat de Catalunya and the Danish
Research Council.}
\address{Departament de Matem\`atiques, Universitat Aut\`onoma de
  Barcelona, 08193, Bellaterra (Barcelona), Spain}
\email{para@mat.uab.es} \email{perera@mat.uab.es}
\address{Department of Mathematics, University of Copenhagen,
  Universitetsparken 5, 2100 Copenhagen \O, Denmark}
\email{gkped@math.ku.dk}
\date{} \dedicatory{} \commby{}
 \subjclass{}
\begin{abstract}
We prove a new extension result for $QB-$rings that allows us to examine
extensions of rings where the ideal is purely infinite and simple. We
then use this result to explore various constructions that provide new
examples of $QB-$rings. More concretely, we show that a surjective
pullback of two $QB-$rings is usually again a $QB-$ring. Specializing to
the case of an extension of a semi-prime ideal $I$ of a unital ring $R$,
the pullback setting leads naturally to the study of rings whose
multiplier rings are $QB-$rings. For a wide class of regular rings, we
give necessary and sufficient conditions for their multiplier rings to be
$QB-$rings. Our analysis is based on the study of extensions and the use
of non-stable $K-$theoretical techniques.
\end{abstract}

\maketitle

\section*{Introduction}
The main theme of the paper \cite{qb} was the study of categorical
properties of the class of $QB-$rings. These constitute a substantial
enlargement of the class of rings with stable rank one (in the sense of
Bass), but enjoy similar nice properties. For example, the property of
being a $QB-$ring passes to matrices and corners and hence is Morita
invariant (\cite[\S 6]{qb}). It also passes to quotients and ideals (with
a suitable formulation of this concept in the non-unital case). Examples
of $QB-$rings include rings whose stable ranks may be infinite: all purely
infinite simple rings (\cite[Proposition 3.10]{qb}), the ring
$\mathrm{End}_D(V)$ where $V$ is a vector space of infinite dimension
over a division ring $D$, the algebra $\mathbb{B}(K)$ of all row- and
column-finite matrices over a field $K$ (\cite[\S 8]{qb}), algebraic
analogues of the Toeplitz algebra (\cite[Examples 7.16 and 7.17]{qb}) and
many more.

The $QB-$notion has its origin in the more geometric concept of extremal
richness for \Cs\ developed by L. G. Brown and the second author in a
series of papers (\cite{bpcrelle}, \cite{bpms}, \cite{bpi} and
\cite{bppre}), where a \C\ $A$ is said to be {\it extremally rich} if the
unit ball of the algebra equals the convex hull of the set of its extreme
points; equivalently, $A$ is extremally rich if the set of so-called
quasi-invertible elements is dense in $A$ (in the norm-topology).

In the algebraic setting, the condition that qualifies a unital ring to
be a $QB-$ring is the following:

For any $a$ and $b$ satisfying $Ra+Rb=R$, there exists an element $y$ in
$R$ such that $a+yb\in R_q^{-1}$.

Here $R_q^{-1}$ is the set of quasi-invertible elements of $R$, which in
the prime case reduces to the union of the left and right invertible
elements (see below for the precise definition). If we substitute
$R_q^{-1}$ by the set $R^{-1}$ of invertible elements in the definition
above, we recover the notion of stable rank one. The precise statement
asserts that a $QB-$ring $R$ has stable rank one if and only if
$R_q^{-1}=R^{-1}$ (see \cite[Proposition 3.9]{qb}). Despite this formal
replacement of the set of invertible elements by the set of
quasi-invertible ones, the work with $QB-$rings requires a larger display
of technology. Observe also that although our definition favours left
multiplication, the concept is left-right symmetric (\cite[Theorem
3.5]{qb}).

For any subset $E$ of a ring $R$ we define $\cl (E)$ to be the set of
elements $a$ in $R$ such that $E\cap (a+Rb)\neq\emptyset$ whenever
$Ra+Rb=R$. Then it follows easily that $R$ is a $QB-$ring (respectively,
$R$ has stable rank one) if and only if $\cl (R_q^{-1})=R$ (respectively,
$\cl (R^{-1})=R$). In this paper (as in \cite{qb}) we shall only be
concerned with applying the operation $\cl$ to the set $R_q^{-1}$;
however, this operation has remarkable similarities to a real closure in a
topological space, especially in the case of \Cs, see \cite{app2}. Indeed,
as we show in \cite[\S 9]{qb} (see also \cite[Proposition 4.6]{app2}),
for a \C\ $A$ we always have $\cl(A_q^{-1})=\ol{A_q^{-1}}$. Hence the
\Cs\ that are $QB-$rings are precisely the extremally rich ones, and thus
a whole new class of examples is available to us.

The behaviour of $QB-$rings under extensions is considerably more
complicated than that of rings with stable rank one. In \cite[Theorem
7.2]{qb}, we give necessary and sufficient conditions for an extension of
$QB-$rings to be a $QB-$ring, which are easily verifiable when one of the
rings in an extension has stable rank one, but in general are of a
technical nature. One of the objectives of this paper is to obtain a
suitable reformulation of the extension result that allows us to use it
in other circumstances as well. This we do in Theorem \ref{e3}, and as a
first application we study extensions of purely infinite simple rings
(Theorem \ref{e8}).

A second application of this extension result is carried out to analyse
examples constructed by means of pullbacks. Hence we obtain (mild)
sufficient conditions ensuring that a pullback of two $QB-$rings is again
a $QB-$ring (Theorem \ref{e10}). Consequently, every (finite) subdirect
product of $QB-$rings is a $QB-$ring.

Since any extension of a (non-unital) semi-prime ideal $I$ of a ring $R$
can be viewed as a pullback of the quotient ring $R/I$ and the multiplier
ring of $I$, we study in the final section of the paper certain classes of
von Neumann regular rings whose multiplier rings are $QB-$rings. In turn,
this is done via a characterisation of the $QB-$property in terms of the
monoid of isomorphism classes of finitely generated projective modules
given in \cite[\S 8]{qb}. Therefore, the techniques developed in
\cite{ap}, \cite{percan} and \cite{perjot} to study multiplier rings of
von Neumann regular rings and \Cs\ with real rank zero are available. We
prove that the multiplier ring $\mulr$ of a simple regular ring $R$ with a
countable unit is $QB$ if and only if $R$ is either purely infinite
simple or artinian. In other cases, we study the presence of the
$QB-$property in the quotient $\mulr/R$, that is reflected on a geometric
condition on the state space of $R$.

\section{Extensions}
\begin{defis}
{\rm Let us start by recalling the concept of quasi-invertibility. We say
that elements $x$ and $y$ in a unital ring $R$ are {\it centrally
orthogonal} provided that $xRy=yRx=0$, and we write $x\perp y$. An
element $u$ of a unital ring $R$ is {\it quasi-invertible}, in symbols
$u\in R_q^{-1}$, provided there exists an element $v$ in $R$ such that
\[
(1-uv)\perp (1-vu)\,.
\]
Necessarily then $u=uvu$, and we may choose $v=vuv$ (so $u$ and $v$ are
regular elements). We say in this situation that $v$ is a {\it
quasi-inverse} for $u$. As discussed in \cite[\S 2]{qb}, this is a
suitable weakening of the notion of one-sided invertibility, because $u\in
R_q^{-1}$ if and only if there exist centrally orthogonal ideals $I$ and
$J$ of $R$ such that $u+I$ is left invertible in $R/I$ and $u+J$ is right
invertible in $R/J$.

A unital ring $R$ is said to be a {\it $QB-$ring} if whenever $Ra+Rb=R$,
there exists an element $y$ such that $a+yb\in R_q^{-1}$. Since this
strongly involves the unit, the concept has to be reformulated to make
sense in the non-unital case. We shall give here a definition in the
setting of a ring $I$ that sits as a two-sided ideal in a unital ring $R$
(see \cite[\S 4]{qb} for a more complete discussion). We say that $I$ is
a $QB-$ring if whenever $xa-x-a+b=0$ for $x$, $a$ and $b$ in $I$, there
exists an element $y$ in $I$ such that $1-(a-yb)\in R_q^{-1}$. We shall
usually say that $I$ is a {\it $QB-$ideal} of the ring $R$.

If $u\in R_q^{-1}$ with quasi-inverse $v$ we shall refer to the centrally
orthogonal elements $p=1-uv$ and $q=1-vu$ as the {\it defect idempotents}
associated with $u$. Since the quasi-inverse is not unique, this is not
strictly accurate, but by \cite[Theorem 2.3]{qb} any other quasi-inverse
will have the form $v'=v+a(1-uv)+(1-vu)b$ for $a$, $b$ in $R$, and the
new defect idempotents
\[
p'=(1-ua)p\,,\qquad\mbox{and}\qquad q'=q(1-bu)\,.
\]
Since $1-uap\in R^{-1}$ with $(1-uap)^{-1}=1+uap$ and
$p'=(1-uap)p(1+uap)$, we see that $p'$ and $p$ are conjugate; and so are
$q'$ and $q$. Conversely, if $e=wpw^{-1}$ for some $w$ in $R^{-1}$, then
$e$ is a defect idempotent associated with $wu$ (with quasi-inverse
$vw^{-1}$). The precise statement is therefore that the defect
idempotents associated with the (equivalence) class of elements
$\{w_1uw_2\mid w_i\in R^{-1}\}$ in $R_q^{-1}$ are exactly the idempotents
conjugate with $p$ and $q$.}
\end{defis}
\begin{lem}
\label{e1} Let $I$ be a $QB-$ideal of a unital ring $R$, and let $e$ be
an idempotent of $R$. Then $eIe$ is a $QB-$ideal of the unital ring $eRe$.
\end{lem}
\begin{proof}
We follow the definition of a non-unital $QB-$ring given above. Assume
that
\[
(exe)(eae)-exe-eae+ebe=0\,,
\]
for some elements $x$, $a$ and $b$ in $I$. Then $(1-exe)(1-eae)+ebe=1$ in
$R$. Since $I$ is a $QB-$ring, we get an element $y$ in $I$ such that
$1-eae+yebe\in R_q^{-1}$. Note now that the element $1-(1-e)yebe$ is
invertible in $R$, whence
\[
1-eae+eyebe=(1-eae+yebe)(1-(1-e)yebe)\in R^{-1}_q\,.
\]
It is then easy to see that $e-eae+eyebe\in (eRe)^{-1}_q$. Therefore
$eIe$ is a $QB-$ring.
\end{proof}
\begin{defis}
{\rm The notion of $QB-$corner has proved central to the theory of
$QB-$rings (e.g. in showing that this class is stable under matrix
formation). We recall here the main definitions and refer the reader to
\cite[\S 5]{qb} for further details.

If $p$ and $q$ are idempotents in a ring $R$ such that $pRq\ne 0$, we say
that an element $x$ in $pRq$ is {\it quasi-invertible} (and write $x\in
(pRq)_q^{-1}$) provided that $(p-xy)\perp (q-yx)$ for some element $y$ in
$qRp$ (again, $x=xyx$ and we may choose $y=yxy$). We define
$\cl^{\sim}((pRq)_q^{-1})$ to be the set of all elements $a$ in $pRq$
such that whenever we have an equation $xa+b=q$ with $x$ in $qRp$ and $b$
in $qRq$, there exists an element $y$ in $pRq$ satisfying that $a+yb\in
(pRq)^{-1}_q$. The set $\ccr^{\sim}((pRq)^{-1}_q)$ is defined
symmetrically. A skew corner $pRq$ is said to be a {\it $QB-$corner}
provided that $\cl^{\sim} ((pRq)^{-1}_q)=pRq$ and
$\ccr^{\sim}((pRq)^{-1}_q)=pRq$. Of course, if $p=q$, then $pRp$ is a
$QB-$corner if and only if it is a $QB-$ring (\cite[Corollary 5.7]{qb}).}
\end{defis}

Recall that two idempotents $e$ and $f$ are said to be (Murray-von
Neumann) {\it equivalent}, in symbols $e\sim f$, provided there exist
elements $x$ in $eRf$ and $y$ in $fRe$ such that $e=xy$ and $f=yx$. We
also write $e\lesssim f$ provided $e\sim f'$ for some idempotent $f'$ such
that $f'\leq f$.

Let $u$ be an element in $R_q^{-1}$ with quasi-inverse $v$, so that
$u=uvu$, $v=vuv$ and if $p=uv$ and $q=vu$, we have $(1-p)\perp (1-q)$. Let
$w$ be an element in $(qRq)_q^{-1}$ with quasi-inverse $w'$, satisfying
the analogous equations $w=ww'w$, $w'=w'ww'$ and $(q-ww')\perp (q-w'w)$.
Observe that $p-uww'v=u(q-ww')v$ is an idempotent equivalent to $q-ww'$,
whence also $(p-uww'v)\perp (q-w'w)$. Set $p_1=1-p$, $p_2=p-uww'v$,
$q_1=q-w'w$ and $q_2=1-q$. We have the following orthogonality relations:
\begin{equation} \label{ai}
\begin{split}
 & p_2\perp q_1\,,\qquad p_1\perp q_2\,,\\
 & p_1p_2=p_2p_1=q_1q_2=q_2q_1=0\,,\\
 & (p_1+p_2)uw=uw(q_1+q_2)=0\,,\\
 & (q_1+q_2)w'v=w'v(p_1+p_2)=0\,.
\end{split}
\tag{$*$}
\end{equation}
Our interest in the next lemma is the study of elements of the form
$uw+t_1+t_2$, with $t_1$ in $p_1Rq_1$ and $t_2$ in $p_2Rq_2$.
\begin{lem}
\label{caguen} Let $R$ be a unital ring and let $a=uw+t_1+t_2$, where
$t_1\in p_1Rq_1$ and $t_2\in p_2Rq_2$, according to the above notation.
\begin{itemize}
\item[(i)]
If $t_1\in (p_1Rq_1)^{-1}_q$ and $t_2\in (p_2Rq_2)^{-1}_q$ then $a\in
R_q^{-1}$\,.
\item[(ii)]
If $t_1\in\cl^{\sim} ((p_1Rq_1)^{-1}_q)$ and $t_2\in\cl^{\sim}
((p_2Rq_2)^{-1}_q)$ then $a\in\cl (R_q^{-1})$. In particular, this holds
if $p_1Rq_1$ and $p_2Rq_2$ are $QB-$corners\,.
\end{itemize}
The assertions {\rm (i)} and {\rm (ii)} also hold if $t_i=0$ and $p_i\perp
q_i$ for $i=1$ or $2$.
\end{lem}
\begin{proof}
(i). Take quasi-inverses $t_1'$ in $q_1Rp_1$ and $t_2'$ in $q_2Rp_2$ for
$t_1$ and $t_2$ respectively. Then, by using the relations \eqref{ai} we
easily check that the element $a'=w'v+t_1'+t_2'$ is a quasi-inverse for
$a$. Indeed, we get that
\[
aa'=uww'v+t_1t_1'+t_2t_2'\qquad\mbox{and}\qquad a'a=w'w+t_1't_1+t_2't_2\,.
\]
Thus
\[
1-aa'=1-p-t_1t_1'+p-uww'v-t_2t_2'=(p_1-t_1t_1')+(p_2-t_2t_2')
\]
and
\[
1-a'a=q-w'w-t_1't_1+1-q-t_2't_2=(q_1-t_1't_1)+(q_2-t_2't_2)\,.
\]
Since $p_1\perp q_2$, $p_2\perp q_1$, $(p_1-t_1t_1')\perp (q_1-t_1't_1)$ and
$(p_2-t_2t_2')\perp (q_2-t_2't_2)$, we conclude that $(1-aa')\perp (1-a'a)$.

(ii). Assume that $xa+b=1$ for some $x$ and $b$ in $R$. Left and right
multiplication by $q_1$, coupled with the fact that $aq_1=t_1$ yields
\[
q_1xt_1+q_1bq_1=q_1\,,
\]
which we can rewrite as $(q_1xp_1)(p_1t_1q_1)+q_1bq_1=q_1$. Since
$t_1\in\cl^{\sim} ((p_1Rq_1)^{-1}_q)$, there exists an element $c_1$ in
$p_1Rq_1$ such that $t_1+c_1bq_1\in (p_1Rq_1)^{-1}_q$. Similarly (using
the hypothesis that $t_2\in\cl^{\sim} ((p_2Rq_2)^{-1}_q)$), we find an
element $c_2$ in $p_2Rq_2$ such that $t_2+c_2bq_2\in (p_2Rq_2)^{-1}_q$.
Recalling that $q_1=q-w'w$ and $q_2=1-q$, we have
\begin{align*}
c_1b & = c_1b(1-q_1)+c_1bq_1\\
& =c_1b(1-q)+c_1bw'w+c_1bq_1\\
& = c_1bq_2+c_1bw'w+c_1bq_1\,.
\end{align*}
Now $c_1\in p_1R$ and by \eqref{ai} we have $p_1\perp q_2$, so $c_1bq_2=0$
and therefore $c_1b=c_1bw'w+c_1bq_1$. Similarly (using this time that
$p_2\perp q_1$), $c_2b=c_2bw'w+c_2bq_2$. Using these relations, we find
that
\begin{align*}
a+(c_1+c_2)b & = uw+t_1+t_2+c_1bw'w+c_1bq_1+c_2bw'w+c_2bq_2\\
& = uw+(c_1+c_2)bw'vuw+(t_1+c_1bq_1)+(t_2+c_2bq_2)\\
& = (1+(c_1+c_2)bw'v)uw+(t_1+c_1bq_1)+(t_2+c_2bq_2)\,.
\end{align*}
Set $d=(c_1+c_2)bw'v$. Since $p_1+p_2=1-uww'v$, we have $d(p_1+p_2)=0$,
and so $d^2=0$ (by \eqref{ai}). Therefore the element $d_1=1+d$ is
invertible (with inverse $d_1^{-1}=1-d$) and it follows from the relations
\eqref{ai} that
$d_1^{-1}((t_1+c_1bq_1)+(t_2+c_2bq_2))=(t_1+c_1bq_1)+(t_2+c_2bq_2)$.
Consequently
\[
d_1^{-1}(a+(c_1+c_2)b)=uw+(t_1+c_1bq_1)+(t_2+c_2bq_2)\,.
\]
Since $t_1+c_1bq_1\in (p_1Rq_1)^{-1}_q$ and $t_2+c_2bq_2\in
(p_2Rq_2)^{-1}_q$, it follows from (i) that $d_1^{-1}(a+(c_1+c_2)b)\in
R_q^{-1}$, whence also $a+(c_1+c_2)b\in R_q^{-1}$ and this shows that
$a\in\cl (R_q^{-1})$, as desired.
\end{proof}
\begin{lem}
\label{e2} Let $I$ be a $QB-$ideal of a unital ring $R$, and assume that
for every pair of defect idempotents $p$ and $q$ with $p$ or $q$ in $I$,
either $p\perp q$ or else $pRq$ (hence also $qRp$) is a $QB-$corner. Then
$I+R^{-1}_q\subset \cl(R^{-1}_q)$.
\end{lem}
\begin{proof}
Let $u$ be an element in $R^{-1}_q$ with quasi-inverse $v$, so that
$u=uvu$, $v=vuv$ and if $p=uv$, $q=vu$, we have $(1-p)\perp(1-q)$.

Now take $t$ in $I$ and assume that there is an equation:
\[
x(u-t)+b=1\,,
\]
for some $x$ and $b$ in $R$. By using that $(1-p)t(1-q)=0$, we can rewrite
this as
\[
xuq(1-vt)q-x(1-p)t-xt(1-q)+b=1\,.
\]
Right and left multiplication by $q$ yields $qxuq(1-vt)q-qx(1-p)tq
+qbq=q$.

Observe that $q(1-vt)q\in q-qIq$. Since $qIq$ is a $QB-$ideal of $qRq$ (by
Lemma \ref{e1}), there exists therefore an element $z$ in $qIq$ (using
\cite[Lemma 4.6]{qb} if necessary) such that the element
\[
w=q(1-vt)q+z(qbq-qx(1-p)tq)\in (qRq)^{-1}_q
\]
(and actually $w\in q+qIq$). By computation (and using again that
$(1-p)t(1-q)=0$)

\begin{align*}
uw & =u-ptq+uzqbq-uzqx(1-p)tq\\
& =u-t+(1-p)tq+pt(1-q)+uzqbq-uzqx(1-p)tq\\
& =u-t+t(1-q)+(1-uzqx(1-p))(1-p)t+uzqbq\\
& =u-t+t(1-q)+w_1^{-1}(1-p)t+uzqb-uzqb(1-q)\\
& =u-t+w_1^{-1}t'(1-q)+w_1^{-1}(1-p)t+uzqb\,,
\end{align*}
where $w_1=1+uzqx(1-p)$ is an invertible element in $R$ and
$t'=w_1(t-uzqb)$. Note that $w_1u=u$. Rearranging the previous equality,
we get
\[
w_1(u-t+uzqb)=uw-t'(1-q)-(1-p)t\,.
\]
Let $a'=uw-t'(1-q)-(1-p)t$, and note that
\[
u-t+uzqb=w_1^{-1}a'\,.
\]
Let $w'$ be a quasi-inverse for $w$ in $qRq$ (and actually $w'\in
q+qIq$). The defect idempotents for $w$ and $w'$ are $q-ww'$ and $q-w'w$,
both belonging to $I$. Note that in fact $q-ww'$ and $q-w'w$ are defect
idempotents of the quasi-invertible element $w+1-q$ (in $R$), with
quasi-inverse $w'+1-q$. Now consider the idempotent $p-uww'v=u(q-ww')v$,
which is equivalent to $q-ww'$, hence also centrally orthogonal to
$q-w'w$. Similarly, we see that it is a defect idempotent for the
quasi-invertible element $1-p+uwv$ (with quasi-inverse $1-p+vw'u$).

Compute that
\begin{align*}
a'&=uw -t'(1-q)-(1-p)t\\
& =uw -pt'(1-q)-(1-p)t\\
 & =uw-uww'vt'(1-q)-(p-uww'v)t'(1-q)-(1-p)t\\
& = uww_2^{-1}-(p-uww'v)t'(1-q)-(1-p)t\,,
\end{align*}
where $w_2=1+w'vt'(1-q)$, as before, is an invertible element. Observe
that $(1-q)w_2=1-q$ and $(1-p)tw_2=(1-p)t$ since $(1-p)\perp (1-q)$.
Therefore,
\[
a'w_2=uw-(p-uww'v)t'(1-q)-(1-p)t\,.
\]
Now,
\begin{align*}
a'w_2  &= uw-(p-uww'v)t'(1-q)-(1-p)tq\\
 &= uw-(p-uww'v)t'(1-q)-(1-p)t(q-w'w)-(1-p)tw'w\\
 &= w_3^{-1}uw-(p-uww'v)t'(1-q)-(1-p)t(q-w'w)\,,
\end{align*}
where $w_3=1+(1-p)tw'v$ is an invertible element. We also have that
$w_3(p-uww'v)=p-uww'v$ and that $w_3(1-p)=1-p$. Thus,
\[
w_3a'w_2=uw-(p-uww'v)t'(1-q)-(1-p)t(q-w'w)\,.
\]
Set $t_1=-(1-p)t(q-w'w)$ and $t_2=-(p-uww'v)t'(1-q)$, and note that the
element $a''=w_3a'w_2=uw+t_1+t_2$ matches with our notation in Lemma
\ref{caguen}. Since $q-w'w$ and $p-uww'v$ belong to $I$, our hypothesis
implies that either $p_1\perp q_1$ (respectively $p_2\perp q_2$) or else
$p_1Rq_1$ (respectively $p_2Rq_2$) is a $QB-$corner, where $p_1$, $q_1$,
$p_2$ and $q_2$ are as in Lemma \ref{caguen}.

It follows from condition (ii) in Lemma \ref{caguen} that $a''\in\cl
(R_q^{-1})$ and so $a'=w_3^{-1}a''w_2^{-1}\in\cl(R_q^{-1})$ and therefore
$u-t+uzqb=w_1^{-1}a'\in\cl(R_q^{-1})$. Since
\[
x(u-t+uzqb)+(1-xuzq)b=1\,,
\]
there is an element $y$ in $R$ such that $u-t+uzqb+y(1-xuzq)b\in
R_q^{-1}$, and so
\[
u-t+(uzq+y(1-xuzq))b\in R_q^{-1}\,.
\]
This shows that $u-t\in\cl(R_q^{-1})$, as desired.
\end{proof}

The immediate benefit of the lemma is the following extension result for
$QB-$rings, analogous to condition (ii) in \cite[Theorem 6.1]{bpcrelle}.
\begin{theor}
\label{e3} Let $R$ be a unital ring and let $I$ be a two-sided ideal of
$R$. Then $R$ is a $QB-$ring if and only if the following conditions are
satisfied:
\begin{itemize}
\item[(i)]
$R/I$ is a $QB-$ring\,;
\item[(ii)]
$(R/I)^{-1}_q=(R^{-1}_q+I)/I$\,;
\item[(iii)]
$I$ is a $QB-$ideal and for every pair of defect idempotents $p$ and $q$
with $p$ or $q$ in $I$, either $p\perp q$ or else $pRq$ is a $QB-$corner.
\end{itemize}
\end{theor}
\begin{proof}
The forward implication holds by combining \cite[Theorem 7.2]{qb} with
\cite[Corollary 4.10]{qb} and \cite[Lemma 6.1]{qb}.

Conversely, condition (iii) above implies that $I+R^{-1}_q\subset \cl
(R^{-1}_q)$, by Lemma \ref{e2}. Thus, by \cite[Theorem 7.2]{qb},
conditions (i)-(iii) imply that $R$ is a $QB-$ring.
\end{proof}
\begin{defis}
{\rm Recall that a unital ring $R$ is an {\it exchange ring} provided
that for any element $x$ in $R$, there exists an idempotent $p$ in $xR$
such that $1-p\in (1-x)R$ (see \cite{gw} and \cite{nic}). Rewriting this
last condition as $1-p=(1-x)(1-s)$ for some $s$ in $R$ we get $p=x+s-xs$,
showing that the role of the unit is superfluous, and this provides a
definition of an exchange ring in the non-unital case (see \cite{aext}).
This class of rings contains many examples, among them all von Neumann
regular rings, all semi-perfect rings and all \Cs\ with real rank zero. It
is a class well-behaved under natural constructions, such as matrix
formation, passage to ideals and quotients, and idempotent-lifting
extensions (see \cite{nic} and \cite{aext}).

A non-zero idempotent $p$ of a ring $R$ is said to be {\it infinite}
provided that there exists an idempotent $q$ in $R$ such that $q<p$ and
$q\sim p$. We say that a simple ring $R$ is {\it purely infinite} if
every non-zero right ideal contains an infinite idempotent. Note that the
formulation chosen here does not require the ring $R$ to be unital. It is
possible to show that the concept is left-right symmetric and that in the
unital case $R$ is purely infinite simple if and only if $R$ is not a
division ring and for any non-zero element $x$ in $R$, there exist
elements $s$ and $t$ in $R$ such that $sxt=1$ (\cite[Theorem 1.6]{agp}).

Purely infinite simple rings satisfy an important technical feature: for
any non-zero idempotents $p$ and $q$ in $R$, we always have $p\lesssim
q$. This was proved in \cite[Proposition 1.5]{agp} and we shall use it
below.

It is known that purely infinite simple \Cs\ have real rank zero, that
is, they are exchange rings (see \cite{zh}, \cite[Theorem 7.2]{agop}).
The analogous result in the purely ring-theoretic setting remains open.}
\end{defis}
\begin{lem}
\label{e4} Let $I$ be a purely infinite simple exchange ring. Assume
furthermore that $I$ is an essential ideal of a unital ring $R$. Then, for
every pair of idempotents $p$ and $q$ in $R$ with $p$ or $q$ in $I$,
either $p\perp q$ or else $pRq$ is a $QB-$corner.
\end{lem}
\begin{proof}
Let $p$ and $q$ be idempotents in $R$, and assume that $p\in I$ and that
$pRq\ne 0$. We have to prove that $\cl^{\sim}((pRq)^{-1}_q)=pRq$ and
$\ccr^{\sim}((pRq)^{-1}_q)=pRq$.

To see this, assume first that we have an equation $xa+b=q$, where $a\in
pRq$, $x\in qRp$ and $b\in qRq$. Since $qIq$ is an exchange ideal of
$qRq$ (see \cite[Proposition 1.3]{aext}), there is an idempotent $r$ in
$(qRq)xa$ such that $q-r\in (qRq)b$. Write $r=txa$ and $q-r=sb$ for some
elements $t$ and $s$ in $qRq$. Note that $r\in I$. If $q-r=0$, then
$q=txa$ and therefore $a\in (pRq)^{-1}_q$, with $tx$ as a quasi-inverse.

Otherwise, $q-r\neq 0$. Since $I$ is essential and simple,
$(q-r)I(q-r)\neq 0$, and since $I$ is purely infinite there exists an
infinite idempotent $e$ in $(q-r)I(q-r)$. It follows that $p-artx\lesssim
e\leq q-r$. Write
\[
p-artx=uv\,,\qquad vu\leq q-r\,,
\]
with $u$ in $(p-artx)R(q-r)$ and $v$ in $(q-r)R(p-artx)$. Now consider the
element $w=u+ar$ ($w\in pRq$), and compute that
\[
w(v+rtx)=(u+ar)(v+rtx)=uv+artx=p\,.
\]
Therefore $w\in (pRq)^{-1}_q$, and in fact
\[
a+(w-a)sb=aq+w(q-r)-a(q-r)=ar+u=w\,.
\]
This shows that $\cl^{\sim}((pRq)^{-1}_q)=pRq$. The proof that
$\ccr^{\sim}((pRq)^{-1}_q)=pRq$ uses similar arguments.
\end{proof}
\begin{prop}
\label{e5} Let $R$ be a unital exchange ring. Then $R$ is a $QB-$ring if
and only if for every pair of idempotents $p$ and $q$ such that $1-p\sim
1-q$, either $p\perp q$ or else $(pRq)^{-1}_q\neq\emptyset$.
\end{prop}
\begin{proof}
The condition is necessary by \cite[Corollary 5.8]{qb}. To establish the
converse, we show that if $x$ is a regular element in $R$, then there
exists a quasi-invertible element $u$ such that $x=xux$ and then apply
\cite[Theorem 8.4]{qb}.

Let $x=xyx$ be a regular element in $R$. Set $p=1-xy$ and $q=1-yx$, and
note that evidently $1-p\sim 1-q$. If $p\perp q$, then $(1-xy)\perp
(1-yx)$, and in this case $x$ and $y\in R^{-1}_q$.

Otherwise, there exists a quasi-invertible element $w$ in $pRq$, by
hypothesis. Now, \cite[Theorem 5.5]{qb} implies that $x+w\in R^{-1}_q$
(with quasi-inverse $y+w'$, where $w'\in qRp$ and is a quasi-inverse for
$w$). Set $u=y+w'$ and compute that $xux=x(y+w')x=xyx=x$.
\end{proof}

In the non-unital case, this characterisation does not seem to hold in
exactly the same way, but still some implications are true:
\begin{lem}
\label{e6} Let $I$ be a two-sided ideal of a unital ring $R$. Consider
the following conditions:
\begin{itemize}
\item[(i)]
For all idempotents $p$ and $q$ in $I$ such that $1-p\sim 1-q$, either
$p\perp q$ or else $pRq$ is a $QB-$corner.
\item[(ii)]
For all idempotents $p$ and $q$ in $I$ such that $1-p\sim 1-q$, either
$p\perp q$ or else $(pRq)^{-1}_q\neq\emptyset$\,.
\item[(iii)]
If $x\in I$ such that $1-x$ is a regular element, then $1-x$ extends to a
quasi-invertible element.
\item[(iv)]
$I$ is a $QB-$ring.
\end{itemize}
We always have {\rm (i)} $\Rightarrow$ {\rm (ii)} $\Rightarrow$ {\rm
(iii)} and {\rm (iv)} $\Rightarrow$ {\rm (iii)}. If $I$ is an exchange
ideal, then also {\rm (iii)} $\Rightarrow$ {\rm (iv)}.
\end{lem}
\begin{proof}
Evidently {\rm (i)} implies {\rm (ii)}. The implication {\rm (iv)}
$\Rightarrow$ {\rm (iii)} follows from \cite[Theorem 4.9]{qb} and
\cite[Theorem 5.10]{qb}.

{\rm (ii)} $\Rightarrow$ {\rm (iii)}. Let $x$ be an element of $I$, and
assume that $1-x$ is regular. Thus we may write $1-x=(1-x)(1-y)(1-x)$ for
some $y$ in $R$ (and actually $y\in I$). Now the argument in the previous
proposition carries through.

The argument that {\rm (iii)} implies {\rm (iv)} if $I$ is an exchange
ideal is very similar to the one used in the implication {\rm (ii)}
$\Rightarrow$ {\rm (i)} in \cite[Theorem 8.4]{qb}, and we therefore omit
the details.
\end{proof}

It was proved in \cite[Proposition 3.10]{qb} that all purely infinite
simple unital rings are $QB-$rings. We have only been able to establish
the parallel result in the non-unital case under the additional assumption
that the ring is exchange.
\begin{corol}
\label{e7} Let $I$ be a purely infinite, simple ideal of a unital ring
$R$. If $I$ is an exchange ideal then $I$ is a $QB-$ideal.
\end{corol}
\begin{proof}
Take idempotents $p$ and $q$ in $I$, and assume that $pRq\neq 0$. Then,
in particular, both $p$ and $q$ are non-zero, so $p\lesssim q$. Thus
$p=xy$ and $yx\leq q$ for some $x$ in $pRq$ and $y$ in $qRp$. Clearly
this implies that $x\in (pRq)^{-1}_q$ (with $y$ as its quasi-inverse).
Now the result follows from Lemma \ref{e6}.
\end{proof}
\begin{theor}
\label{e8} Let $I$ be a purely infinite, simple and essential ideal of a
unital ring $R$. If $I$ is an exchange ring, then $R$ is a $QB-$ring if
and only if $R/I$ is a $QB-$ring and $(R/I)^{-1}_q=(R/I)^{-1}_r\cup
(R/I)^{-1}_l$.
\end{theor}
\begin{proof}
Assume first that $R$ is a $QB-$ring. Then $R/I$ is a $QB-$ring by
\cite[Corollary 3.8]{qb}. Moreover, since $I$ is simple and essential,
$R$ is prime and so we have $R^{-1}_q=R^{-1}_l\cup R^{-1}_r$. Now, if
$u\in (R/I)^{-1}_q$, we first lift $u$ to a quasi-invertible $x$ in $R$
(by \cite[Theorem 7.2]{qb}), which is in fact left or right invertible,
as we have just observed. Therefore $u$, being the image of $x$ is also
left or right invertible.

Conversely, assume that $R/I$ is a $QB-$ring and that
$(R/I)^{-1}_q=(R/I)^{-1}_r\cup (R/I)^{-1}_l$. By Corollary \ref{e7}, $I$
is a $QB-$ideal and by Lemma \ref{e4}, for any pair of (defect)
idempotents $p$ and $q$ with $p$ or $q$ in $I$, either $p\perp q$ or else
$pRq$ is a $QB-$corner.

Hence, in order to apply Theorem \ref{e3} we only have to show that
one-sided invertible elements lift. For any $x$ in $R$, denote by $\ol x$
its equivalence class in $R/I$. Assume that $\ol{xy}=1$. Since $I$ is an
exchange ideal of $R$, we may use \cite[Lemma 2.1]{aext} to find elements
$a$ and $b$ in $R$ such that $a=aba$, $b=bab$ and $\ol a =\ol x$, $\ol
b=\ol y$. Without loss of generality, we may assume that the idempotents
$1-ab$ and $1-ba$ are both non-zero. Since $(1-ba)I(1-ba)\neq 0$ and $I$
is purely infinite and simple, there is an infinite idempotent $r$ such
that $1-ab\lesssim r\leq 1-ba$. Thus, if we write $1-ab=ts$ and $st\leq
1-ba$ for some elements $s$ and $t$ in $I$, we get that $t+a$ is right
invertible (with $s+b$ as a right inverse), and is the required lift for
$a$.
\end{proof}
\begin{corol}
\label{unomas} Let $R$ be a unital ring and let $I$ be an ideal of $R$
satisfying the hypotheses in Theorem \ref{e8}. If moreover $R/I$ is a
prime ring, then $R$ is a $QB-$ring if and only if $R/I$ is a $QB-$ring.
\end{corol}

\section{PullBacks}

\begin{defis}
\label{yata} {\rm For every unital homomorphism $\pi:R\fl S$ we obviously
have $\pi(R^{-1})\subset S^{-1}$, and moreover $\pi(u)^{-1}=\pi(u^{-1})$
for every $u$ in $R^{-1}$. By contrast, we need not have
$\pi(R_q^{-1})\subset S_q^{-1}$ (when $\pi$ is not surjective). For a
simple counterexample, let $u$ denote the unilateral shift on the Hilbert
space $l^2$ and take $\pi$ to be the diagonal embedding of
$\mathbb{B}(l^2)\oplus\mathbb{B}(l^2)$ into
$M_2(\mathbb{B}(l^2))=\mathbb{B}(l^2\oplus l^2)$. The element $u\oplus
u^*$ is invertible in $\mathbb{B}(l^2)\oplus\mathbb{B}(l^2)$, but its
image in $\mathbb{B}(l^2\oplus l^2)$ is neither an isometry nor a
co-isometry, and therefore not quasi-invertible.

In \cite[6.7]{qb} we defined a subring $R$ of $S$ to be {\it primely
embedded} if $p\perp q$ in $R$ implies $p\perp q$ in $S$ for any pair of
idempotents $p$, $q$ in $R$. Evidently this implies that $R_q^{-1}\subset
S_q^{-1}$ if both $R$ and $S$ are unital with the same unit, and in the
non-unital case we get $R_q^\circ\subset S_q^\circ$ for the
quasi-adversible elements, cf. \cite[4.1]{qb}.

We now define a subring $R$ of $S$ to be {\it quasi-primely embedded} if
$R_q^{-1}\subset S_q^{-1}$, when both $R$ and $S$ are unital with the
same unit. If $R$ is non-unital, but $S$ is unital, we demand that
$R_q^\circ\subset S_q^{-1}$; and if both are non-unital we demand that
$R_q^\circ\subset S_q^\circ$. For $C^*-$algebras this condition was
termed {\it extreme point preservation} in \cite[5.1]{bpcrelle}, but this
name makes no sense in the context of ring theory.

For most applications the concept of prime embeddings will suffice, and
the condition is easy to check; however, in some cases the more
specialized quasi-prime embeddings are needed. Note for example that the
ring $K\oplus K$ (where $K$ is any field) embeds as the diagonal in
$M_2(K)$, but the two centrally orthogonal idempotents $(1,0)$
and $(0,1)$ in $K\oplus K$ are not centrally orthogonal
in $M_2(K)$. The embedding is, however, quasi-prime since
$(K\oplus K)^{-1}_q=(K\oplus K)^{-1}\subset
M_2(K)^{-1}$.

It should be noted that if $R$ is quasi-primely embedded in $S$ (and both
are unital) and $u\in R_q^{-1}$, then any quasi-inverse $v$ for $u$ in
$R$ will also be a quasi-inverse for $u$ in $S$. This follows from
\cite[Theorem 2.3]{qb}, noting that $u\in S_q^{-1}$ by assumption and
that $v$ is, after all, a partial inverse for $u$ in $R$, hence also in
$S$ (i.e. $uvu=u$).}
\end{defis}
Our next result is an application of the same Theorem 2.3 to the effect
that if one can lift a quasi-invertible element from a quotient ring,
then one can also lift any of its quasi-inverses.

\begin{prop}
\label{e9} If $I$ is an ideal in a unital $QB-$ring $R$ and $x$, $y$ is a
pair of elements in $(R/I)^{-1}_q$ that are quasi-inverses for each
other, there is a lift of $x$, $y$ to a pair of elements in $R_q^{-1}$
which are still quasi-inverses to each other.
\end{prop}
\begin{proof}
Let $\pi:R\to R/I$ denote the quotient morphism. By \cite[Proposition
7.1]{qb} we can find an element $v$ in $R_q^{-1}$ such that $\pi(v)=y$.
Now choose a quasi-inverse $u$ for $v$ and an element $a$ in $R$ with
$\pi(a)=x$. By \cite[Theorem 2.3]{qb} the element \[
w=u+a(1-vu)+(1-uv)a\in R_q^{-1}\] and $w$ is a quasi-inverse for $v$.

Evidently, \[\pi(v)\pi(a)\pi(v)=yxy=y=\pi(v)\,.\] But we also have
$\pi(v)\pi(u)\pi(v)=\pi(vuv)=\pi(v)$.

Thus, both $x$ ($=\pi(a)$) and $\pi(u)$ are partial inverses to $y$
($=\pi(v)$), whence, by \cite[Theorem
2.3]{qb}\[x=\pi(u)+x(1-y\pi(u))+(1-\pi(u)y)x=\pi(u+a(1-vu)+(1-uv)a)=\pi(w)\,.\]
It follows that $w$, $v$ is the required lift of the pair $x$, $y$.
\end{proof}

Given a commutative diagram of rings with connecting ring homomorphisms:
\[\begin{CD} R @>{\gamma}>> B\\
@V{\delta}VV  @VV{\beta}V\\
A @>{\alpha}>> C
\end{CD}\]
there is a universal solution, viz. the pullback ring
\[A\oplus_C B=\{(a,b)\in A\oplus B\mid \alpha(a)=\beta(b)\}\]
such that there exists a unique morphism $\sigma:R\to A\oplus_C B$ with
$\delta=\pi_1\circ\sigma$ and $\gamma=\pi_2\circ\sigma$, where $\pi_1$,
$\pi_2$ denote the projections on the coordinates of $A\oplus_C B$.
Indeed, $\sigma(x)=(\delta(x), \gamma(x))$.

It is easy to verify, cf. the proof of \cite[Proposition 3.1]{ppp}, that
a commutative diagram as above is a pullback if and only if the following
conditions hold:
\begin{itemize}
\item[(i)]
$\ker\gamma\cap\ker\delta=\{0\}$
\item[(ii)]
$\beta^{-1}(\alpha(A))=\gamma(R)$
\item[(iii)]
$\delta(\ker\gamma)=\ker\alpha$
\end{itemize}
We shall be particularly interested in pullback diagrams in which one of
the morphisms, say $\alpha$, is surjective. By condition (ii) this implies
that also $\gamma$ is surjective, and thus, by (i) and (iii), $\delta$ is
an isomorphism between $\ker\gamma$ and $\ker\alpha$. Setting
$I=\ker\gamma=\ker\alpha$ this means that such a pullback diagram is
described by two extensions with a common ideal, viz.
\[\begin{CD}
0 @>>> I @>>>  R @>{\gamma}>> B @>>> 0\\
@. @| @V{\delta}VV @VV{\beta}V @.\\
0 @>>> I @>>> A @>{\alpha}>> C @>>> 0
\end{CD}\]
Conversely, by conditions (i)-(iii), each such commutative diagram of
extensions describes a pullback.

It is perhaps worth mentioning that any subdirect product $R$ of rings
$A$ and $B$ can be described by a surjective pullback diagram as above.
Indeed, if $\pi_1$ and $\pi_2$ denote the coordinate projections of $R$,
let
\[I=\ker\pi_2=\{(x,0)\in R\mid x\in A\}\,.\]
Evidently the projection $\pi_1:R\to A$ is injective on $I$. To verify
that $\pi_1(I)$ is an ideal in $A$ take any $a$ in $A$ and $x$ in
$\pi_1(I)$. Since $\pi_1$ is surjective, $(a,b)\in R$ for some $b$ in
$B$, whence $(a,b)(x,0)=(ax,0)\in \ker\pi_2$, so $ax\in\pi_1(I)$ (and
similarly $xa\in\pi_1(I)$). Identifying $I$ and $\pi_1(I)$ and putting
$C=A/I$ we have a surjective pullback diagram. The morphism $\beta:B\to
C$ is defined on equivalence classes as
\[\beta((a,b)+I)=a+I\,,\quad (a,b)\in R\,.\]
Note that we only need the surjectivity of $\pi_2$ in order to define
$\beta$ on all of $B$, not just on $\pi_2(R)$. Note also that if
$J=\ker(\pi_1)$($=\ker\delta$), then $I\cap J=\{0\}$ in $R$ and
$C=A/I=B/J$.
\begin{theor}
\label{e10} Consider a surjective pullback diagram of unital rings and
morphisms
\[\begin{CD}
0 @>>> I @>>>  R @>{\gamma}>> B @>>> 0\\
@. @| @V{\delta}VV @VV{\beta}V @.\\
0 @>>> I @>>> A @>{\alpha}>> C @>>> 0
\end{CD}\]
where both rows are extensions. If $A$ and $B$ are $QB-$rings and $\beta
(B)$ is (quasi-) primely embedded in $C$, then $R$ is a $QB-$ring and
$\delta(R)$ is (quasi-) primely embedded in $A$.
\end{theor}
\begin{proof}
If $p$ and $q$ are centrally orthogonal idempotents in $R$, then
$\gamma(p)\perp\gamma(q)$ in $B$, since $\gamma$ is surjective. Since
$\beta(B)$ is primely embedded in $C$ this implies that
$\beta(\gamma(p))\perp\beta(\gamma(q))$ in $C$. Thus
\[\alpha(\delta(p)A\delta(q))=\beta(\gamma(p))C\beta(\gamma(q))=0\,,\]
so $\delta(p)A\delta(q)\subset I$. But then
\[\delta(p)A\delta(q)=\delta(p)\delta(p)A\delta(q)\delta(q)\subset\delta(p)I\delta(q)=\delta(pIq)=0\,.\]
Similarly $\delta(q)A\delta(p)=0$, so $\delta(p)\perp\delta(q)$ in $A$,
proving that $\delta(R)$ is primely embedded in $A$.

If $\beta(B)$ is quasi-primely embedded in $C$ and $u\in R_q^{-1}$ with
quasi-inverse $v$, then $\gamma (u)\in B_q^{-1}$ since $\gamma$ is
surjective, and thus $\beta(\gamma (u))\in C_q^{-1}$ by assumption.
Moreover, as noted in \ref{yata}, the element $\beta(\gamma(v))$ is a
quasi-inverse for $\beta(\gamma(u))$. It follows that if $p=1-uv$ and
$q=1-vu$ then
\[
\alpha(\delta(p)A\delta (q))=\beta(\gamma(p))C\beta(\gamma(q))=0\,,
\]
so that $\delta(p)A\delta(q)\subset I$. However,
\[
\delta(p)I\delta(q)=pIq=0
\]
since $\delta|I$ is an isomorphism and $p\perp q$ in $R$. Consequently
$\delta(p)\perp \delta (q)$ in $A$, whence $\delta(u)\in A_q^{-1}$, as
desired.

To show that $R$ is a $QB-$ring we shall verify the conditions (i)-(iii)
in Theorem \ref{e3} for the extension $0\to I\to R\to B\to 0$. By
assumption $B$ is a $QB-$ring, and so is $I$, being (isomorphic to) an
ideal in the $QB-$ring $A$, cf. \cite[Corollary 4.10]{qb}. To verify that
quasi-invertible elements lift, consider $u$ in $B_q^{-1}$ and choose a
quasi-inverse $v$ for $u$. Then $(1_B-uv)\perp(1_B-vu)$ in $B$, whence
$(1_C-\beta(uv))\perp(1_C-\beta(vu))$ in $C$, since $\beta(B)$ is
(quasi-) primely embedded in $C$. Consequently $\beta(u)\in C_q^{-1}$ with
quasi-inverse $\beta(v)$. Since $A$ is a $QB-$ring we can use Proposition
\ref{e9} to lift the pair $\beta(u)$, $\beta (v)$ to a pair $x$, $y$ in
$A_q^{-1}$ which are quasi-inverses for each other. By construction
$(u,x)\in R$ and $(v,y)\in R$. Moreover,
\[(1_R-(u,x)(v,y))\perp (1_R-(v,y)(u,x))\]
in $R$ (indeed in $A\oplus B$), so $(u,x)\in R_q^{-1}$, as desired.

Finally, consider two defect idempotents $p$ and $q$ in $R$ and assume
that $p\in I$. Choose $u$ and $x$ in $R_q^{-1}$ with quasi-inverses $v$
and $y$, respectively, such that
\[p=1-uv\mbox{ and }q=1-yx\,.\]
To prove that $pRq$ is a $QB-$corner (if $p$ and $q$ are not centrally
orthogonal) consider an equation $xa+b=q$ with $x$ in $qRp$, $a$ in $pRq$
and $b$ in $qRq$. Applying $\delta$ we obtain the equation
$\delta(x)\delta(a)+\delta(b)=\delta(q)$ and since $\delta(p)$ and
$\delta (q)$ are defect idempotents in $A$, because $\delta(R)$ is
(quasi-) primely embedded in $A$, we know from \cite[Lemma 6.1]{qb} that
$\delta(p)A\delta(q)$ is a $QB-$corner in $A$. There is therefore an
element $y$ in $\delta(p)A\delta(q)$ such that
\[\delta(a)+y\delta(b)\in (\delta(p)A\delta(q))^{-1}_q\,.\]
Now observe that
\[\delta(p)A\delta(q)\subset\delta(p)I\delta(q)=\delta(pIq)\subset\delta(pRq)\,,\]
so that $y=\delta(z)$ for some element $z$ in $pRq$. The same computation
shows that $\delta$ is an isomorphism between $pRq$ and
$\delta(p)A\delta(q)$. Consequently, the statement
\[\delta(a+zb)\in (\delta(pRq))^{-1}_q\]
is equivalent to the assertion that
\[a+zb\in (pRq)^{-1}_q\,.\]
We have proved that $pRq\subset\cl^{\sim}((pRq)^{-1}_q)$ and a symmetric
argument shows that also $pRq\subset\ccr^{\sim}((pRq)^{-1}_q)$, whence
$pRq$ is a $QB-$corner in $R$, as desired.
\end{proof}
\begin{corol}
\label{e11} Any finite subdirect product of $QB-$rings is again a
$QB-$ring.
\end{corol}

\section{Multiplier rings}

Further applications of Theorem \ref{e10} require the use of multiplier
rings. Let $I$ be a non-unital semi-prime ring. Recall that the
multiplier ring $\mathcal{M}(I)$ of $I$ is a unital ring that contains
$I$ as a two-sided ideal, and is universal with this property (that is,
$\mathcal{M}(I)$ is the unique solution to the universal problem of
adjoining a unit to $I$), see \cite{ho}, \cite{ap}. In $C^*-$algebra
theory, this ring has been studied and used in many instances, see, e.g.
\cite{percan} and the references therein.

The connection between our earlier results and multiplier rings is the
following. If $I$ is a semi-prime ideal of a unital ring $R$, then the
universal property of the multiplier ring $\mathcal{M}(I)$ combined with
a standard diagram chase argument ensures that we have a commutative
diagram with exact rows
\begin{equation}
\begin{CD}
0 @>>> I @>>>  R @>{\gamma}>> R/I @>>> 0\\
@. @| @V{\delta_{I,R}}VV @VV{\beta_{I,R}}V @.\\
0 @>>> I @>>> \mathcal{M}(I) @>{\alpha}>> \mathcal{M}(I)/I @>>> 0
\end{CD}
\tag{$\dagger$}
\end{equation}
where the maps $\delta_{I,R}$, $\beta_{I,R}$ are determined by $I$ and
$R$, and the right hand square is a pullback. For simplicity, we write
$\delta=\delta_{I,R}$ and $\beta=\beta_{I,R}$ in the following. From
Theorem \ref{e10} we immediately obtain:
\begin{prop}
\label{e12} Let $I$ be a semi-prime ideal of a unital ring $R$. If $R/I$
and $\mathcal{M}(I)$ are $QB-$rings, and $\beta(R/I)$ is (quasi-) primely
embedded in $\mathcal{M}(I)/I$, then $R$ is a $QB-$ring.
\end{prop}
By using Corollary \ref{e11} we can reduce our study of the $QB-$property
in an extension to the case where the ideal is essential. To illustrate
this, we remove the essentiality hypothesis in Corollary \ref{unomas} as
follows:
\begin{prop}
\label{e13} Let $I$ be a purely infinite simple ideal of a unital ring
$R$. If $I$ is an exchange ring and $R/I$ is a prime $QB-$ring, then $R$
is a $QB-$ring.
\end{prop}
\begin{proof}
Consider the pullback diagram $(\dagger)$ obtained as in the previous
discussion from the extension
\[
0\rightarrow I\rightarrow R\rightarrow R/I\rightarrow 0\,.
\]
Let $J=\ker\delta$, and note that $I\cap J=0$. Thus the map
$\overline{\delta}:R/J\rightarrow \mathcal{M}(I)$ is injective and $I$ is
isomorphic to $(I+J)/J$. It follows that in the extension
\[
0\rightarrow I\rightarrow R/J\rightarrow R/(I+J)\rightarrow 0\,,
\]
$I$ is essential in $R/J$. Observe that $R/(I+J)$ is a $QB-$ring, since
it is a quotient of $R/I$. If $u\in (R/(I+J))^{-1}_q$, then we can lift
$u$ to a quasi-invertible element $x$ in $R/I$, because $R/I$ is a
$QB-$ring. Since moreover it is a prime ring, $x$ must be either left or
right invertible and hence the same holds for $u$. Therefore $R/J$ is a
$QB-$ring, by Theorem \ref{e8}.

Since now $R$ is a subdirect product of $R/J$ and $R/I$, we conclude by
Corollary \ref{e11} that $R$ is a $QB-$ring.
\end{proof}

In order to use Proposition \ref{e12} effectively, we need to produce
examples of semi-prime rings whose multiplier rings are $QB-$rings. In
$C^*-$algebra theory, this problem has been addressed in some instances,
mainly \cite{lar}, \cite{perjot} and \cite{bppre}. As noted there, and as
we shall see below, this is a quite strong property even for reasonably
well-behaved classes of rings. However, the quotient ring $\mulr/R$ has
far better chances of being a $QB-$ring. In fact, we provide a fairly
simple way of testing whether such a quotient is a $QB-$ring, thus adding
more (constructible) examples to this class.

In the rest of this section, we shall analyse a wide class of non-unital
von Neumann regular rings, and we shall make extensive use of the
techniques developed in \cite{ap}, \cite{qb} and \cite{percan}. By this
approach it is possible to give new and shorter proofs for some of the
results in \cite{lar} and \cite{perjot}.
\begin{defis}
{\rm For any ring $R$ let $M_{\infty}(R)$ be the directed union of
$M_n(R)$ ($n\in\mathbb N$) and denote by $V(R)$ the monoid of equivalence
classes (in the sense of Murray and von Neumann) of idempotents in
$M_{\infty}(R)$, where the addition is defined by direct sums of
representatives. (The equivalence class of an idempotent $p$ in
$M_{\infty}(R)$ is usually denoted by $[p]$.) This monoid can be
preordered with the so-called {\it algebraic pre-ordering}: if $x$, $y\in
V(R)$, we write $x\leq y$ provided there exists $z$ in $V(R)$ such that
$x+z=y$. An {\it order-ideal} of $V(R)$ is a submonoid $S$ of $V(R)$ such
that whenever $x$ and $y$ in $V(R)$ satisfy that $x\leq y$ and $y\in S$,
then $x\in S$. For example, if $I$ is an ideal of $R$, then it is easy to
see that $V(I)$ is an order-ideal of $V(R)$. Under certain assumptions,
the lattice of order ideals of $V(R)$ parametrizes the lattice of
two-sided ideals of the ring (see, e.g. \cite[Theorem 2.7]{ap}).

If $I$ and $J$ are two ideals of $R$, then we have $V(I)\cap V(J)=V(I\cap
J)=V(IJ)$. Given a ring $R$ and an ideal $I$, define a congruence
relation in $V(R)$ by setting $x\sim y$ if and only if $x+z=y+w$, for
some $z$, $w$ in $V(I)$. Denote by $V(R)/V(I)$ the quotient set modulo
this congruence relation, and note that we can equip it with a monoid
structure by setting $\ol{x}+\ol{y}=\ol{x+y}$, where $\ol{x}$ denotes the
congruence class of $x$. It makes perfect sense to try to relate $V(R/I)$
and $V(R)/V(I)$: there is always a natural monoid morphism, induced by
the projection $\pi:R\rightarrow R/I$, which turns out to be an
isomorphism if $R$ is an exchange ring (\cite[Proposition 1.4]{agop}).}
\end{defis}
We note in passing that if $R$ is a unital and semi-prime ring, then
$M_{\infty}(R)$ is always non-unital, and its multiplier ring is the ring
$\mathbb{B}(R)$ of all row- and column-finite matrices over $R$ (see
\cite[Proposition 1.1]{ap}). This is an important and useful example to
have in mind.

Recall that a semi-prime ring $R$ has a {\it countable unit} provided
that there exists an increasing sequence of idempotents $(p_n)$ in $R$
such that $R=\bigcup p_nRp_n$ (see \cite{ap}). If $R$ is a simple regular
ring, it is a longstanding open question to decide whether $R$ must be
purely infinite simple or have stable rank one (\cite{vnrr},
\cite{agop}). It is nevertheless sensible to approach this class by
distinguishing the two ``extremes" of the scale. We first deal with the
purely infinite case. A result that we shall make use of tacitly in the
sequel is the remarkable theorem proved by K. C. O'Meara in \cite{om},
which asserts that if $R$ is a regular ring with a countable unit, then
$\mathcal{M}(R)$ is an exchange ring.
\begin{theor}
Let $R$ be a purely infinite simple regular ring with a countable unit.
Then $\mulr$ is a $QB-$ring.
\end{theor}
\begin{proof}
By \cite[Proposition 2.13]{ap}, we have that $V(\mulr)\cong
V(R)\sqcup\{\infty\}$, where $x+\infty=\infty$ for all $x$ in $V(R)$ and
$\infty+\infty=\infty$. Now, since we have an isomorphism
$V(\mulr)/V(R)\cong \{0,\infty\}$, we see that the only order-ideals in
$V(\mulr)$ are the trivial ones and $V(R)$. Using e.g. \cite[Theorem
2.7]{ap}, we conclude that the only non-trivial ideal of $\mulr$ is $R$.
Moreover, we know that $V(\mulr/R)\cong V(\mulr)/V(R)$; this, coupled
with the fact that $\mulr$ is an exchange ring (\cite[Theorem 2]{om})
allows us to conclude that $\mulr/R$ is a purely infinite simple ring,
hence a $QB-$ring by \cite[Proposition 3.10]{qb}. Finally, the conclusion
follows by applying Theorem \ref{e8} to the extension
\[0\rightarrow R\rightarrow\mulr\rightarrow\mulr/R\rightarrow 0\,.\]
\end{proof}
If, instead of being purely infinite simple, the regular ring $R$ has
stable rank one, then the intricacy of the corresponding monoid
$V(\mathcal{M}(R))$ is much higher. We will impose a further condition on
the ring, namely that $V(R)$ is a {\it strictly unperforated} monoid.
This by definition means that, whenever $nx+z=ny$, for $x$, $y$, $z$ in
$V(R)$ ($z\neq 0$) and $n$ in $\mathbb{N}$, then $x+w=y$ for some
non-zero $w$ in $V(R)$. We remark that no examples are known of simple
regular rings with stable rank one that fail to be strictly unperforated.

Recall that a {\it state} on a monoid $M$ with order-unit $u$ is a monoid
morphism $s:M\rightarrow \mathbb{R}^+$ such that $s(u)=1$. We denote the
set of states by $\mathrm{St}(M,u)$ or by $S_u$ if no confusion may
arise, and note that it coincides with the set of states on $(G(M),u)$,
the Grothendieck group of $M$. The extreme boundary of $S_u$ (as a convex
set) will be denoted by $\partial_e S_u$. Denoting by $\mathrm{Aff}(S_u)$
the Banach space of affine continuous (real-valued) functions defined on
$S_u$, there is always a natural map $\phi_u:M\rightarrow
\mathrm{Aff}(S_u)$, given by evaluation. We define
$\mathrm{LAff}_{\sigma}(S_u)^{++}$ to be the semi-group whose elements
are those affine and lower semi-continuous functions with values on
$\mathbb{R}^{++}\cup\{+\infty\}$ that can be expressed as countable
(pointwise) suprema of an increasing sequence of (strictly positive)
elements from $\mathrm{Aff}(S_u)$. We shall designate by $\mathcal{N}$
the class of simple regular rings $R$ with a countable unit, stable rank
one, strict unperforation and such that the state space $S_u$ is
metrizable for some, hence any, non-zero element $u$ in $V(R)$ (see
\cite{ap}).

We say that a simple ring $R$ is {\it elementary} provided $R$ has
minimal idempotents. Clearly, $R$ is elementary if and only if is
artinian. In fact, if $R$ is elementary and $p$ is a minimal idempotent
in $R$, then $D=pRp$ is a division ring. Set $V=pR$, $W=Rp$. Then $V$
(respectively $W$) is a left (respectively right) $D-$vector space.
Moreover, $(V,W)$ is a pair of dual spaces, with $V\times W\rightarrow D$
given by $\langle px,yp\rangle=pxyp$. One can then show that
$R\cong\mathfrak{F}_W(V)$, the set of finite rank adjointable operators,
and $\mathcal{M}(R)\cong \mathcal{L}_W(V)$, the set of all adjointable
operators (see e.g. \cite[Proposici\'{o} 3.1.4]{pertesi}). In case $R$ has a
countable unit we moreover have that $R\cong M_{\infty}(D)$ and
$\mulr\cong\mathbb{B}(D)$, and this is a $QB-$ring by the arguments in
\cite[\S 8]{qb}.

If $R$ is a ring in the class $\mathcal{N}$ with a countable unit
$(p_n)$, and if $u$ is a non-zero element in $V(R)$, define
$d=\sup\phi_u([p_n])$, and
\[W_{\sigma}^d(S_u)=\{f\in\mathrm{LAff}_{\sigma}(S_u)^{++}\mid f+g=nd \mbox{
for some }g\in\mathrm{LAff}_{\sigma}(S_u)^{++}\mbox{ and }n\in \mathbb{N}
\}\,.\] Now consider the set $V(R)\sqcup W_{\sigma}^d(S_u)$ (where
$\sqcup$ stands for disjoint union), and equip it with a monoid structure
by extending the natural operations in $V(R)$ and $W_{\sigma}^d(S_u)$,
and by defining $x+f=\phi_u(x)+f$, for any $x$ in $V(R)$ and $f$ in
$W_{\sigma}^d(S_u)$. The result that we highlight below for future
reference is the key for understanding the ideal structure of multiplier
rings of rings in the class $\mathcal{N}$ (for a proof, see \cite[Theorem
2.11 ]{ap}, and also \cite[Theorem 2.6]{percan}).
\begin{theor}
\label{plasta} Let $R$ be a ring in the class $\mathcal{N}$, and let $u$
be a non-zero element in $V(R)$. If $R$ is non-elementary, then there is a
monoid isomorphism
\[\varphi:V(\mathcal{M}(R))\rightarrow V(R)\sqcup W_{\sigma}^d(S_u)\,,\]
such that $\varphi([p])=[p]$ if $p\in R$, and
$\varphi([p])=\sup\{\phi_u([q])\mid [q]\in V(R) \mbox{ and } q\lesssim
p\}$ if $p\in\mathcal{M}(R)\setminus R$.
\end{theor}
We refer to the lower semi-continuous, affine function $d$ in Theorem
\ref{plasta} as the {\it scale of $R$}. The ring $R$ is said to have {\it
continuous scale} (resp. {\it finite scale}) if $d$ is continuous (resp.
the restriction of $d$ to the extreme boundary of $S_u$ is finite). We
say that a state $s$ in $S_u$ is {\it infinite} if $d(s)=\infty$. Theorem
\ref{plasta} (and its $C^*-$counterpart) was used in \cite{ap} and
\cite{percan} to study the ideal lattice of multiplier rings for a wide
class of rings. The only item of this analysis we shall be using here is
the existence of the smallest ideal $L(R)$ of $\mulr$ properly containing
$R$, whose monoid of equivalence classes of idempotents is isomorphic
(through the isomorphism of Theorem \ref{plasta}) to
$V(R)\sqcup\mathrm{Aff}(S_u)^{++}$ (see \cite[\S 3]{ap} and
\cite[Proposition 4.1]{percan}).

To appreciate the following results, we note that Theorem \ref{plasta},
coupled with the fact that $\mathcal{M}(R)$ is an exchange ring if $R$ is
regular with a countable unit (\cite{om}), allows to compute the stable
rank of these multiplier rings, with a proof that follows verbatim the
arguments in \cite[Theorem 7.5]{percan} and that, for this reason, we omit
in the theorem below. Therefore, the question of analysing the
``$QB-$ness'' of $\mathcal{M}(R)$ is quite pertinent.
\begin{theor}
\label{plasta2} Let $R$ be a ring in the class $\mathcal{N}$. Let $u$ be a
non-zero element in $V(R)$. Then $\mathcal{M}(R)$ has stable rank $2$ if
and only if $R$ has finite, but not continuous scale; and it is infinite
otherwise.
\end{theor}
\begin{rema}
\label{otia} {\rm Another ingredient that we shall be using is the
characterisation of being a $QB-$ring given in \cite[Theorem 8.7]{qb} and
restated here for the convenience of the reader:

If $R$ is a semi-primitive exchange ring, then $R$ is a $QB-$ring if and
only if the following condition is met in $V(R)$: if $n\cdot
[1]+b_1=n\cdot [1]+b_2$ in $V(R)$ for some $n\geq 1$, then there exist
order-ideals $S_1$, $S_2$ in $V(R)$ such that $S_1\cap S_2=0$, and
elements $x'$, $c_1$, $c_2$ such that $c_i\in S_i$ for $i=1$, $2$, and
$b_1+c_1=b_2+c_2$, and $x'+c_1=n\cdot [1]=x'+c_2$.

Observe that this characterisation has some content only when both $b_1$
and $b_2$ are non-zero. If $b_1=0$, say, we may choose $S_1=V(R)$,
$S_2=0$, and then necessarily $c_2=0$, $c_1=b_2$ and $x'=n\cdot [1]$.}
\end{rema}
In the following results, we shall use the notation $\{s\}'$, for an
extreme point $s$ of $S_u$, to designate the {\it complementary face} of
the face $\{s\}$. This is, by definition, the union of faces in $S_u$
disjoint from $\{s\}$ (see \cite{poag} and \cite{vnrr}).
\begin{theor}
\label{mult1} Let $R$ be a non-elementary ring in the class $\mathcal{N}$,
and let $p$ be a non-zero idempotent in $R$. Then:
\begin{itemize}
\item[(i)]
If $R$ has finite scale, then $\mulr/R$ is a $QB-$ring.
\item[(ii)]
If $R$ has at least two infinite extremal states in $S_{[p]}$, then
$\mulr/R$ is not a $QB-$ring.
\end{itemize}
\end{theor}
\begin{proof}
(i). We have to show that the monoid $V(\mulr/R)$ satisfies the condition
in Remark \ref{otia}. First, we identify this monoid with
$V(\mulr)/V(R)$, and $V(\mulr)$ with $V(R)\sqcup W_{\sigma}^d(S_u)$,
where $u=[p]$ is non-zero in $V(R)$ and $d=\sup\phi_u([p_n])$, where
$(p_n)$ is a countable unit for $R$, as in Theorem \ref{plasta}. Assume
then that $n\ol{d}+\ol{b_1}=n\ol{d}+\ol{b_2}$, for some $n\geq 1$.

By Remark \ref{otia}, we can assume that $\ol{b_i}\neq 0$ for $i=1$, $2$,
that is, $b_i\notin V(R)$ and so there exist elements $x$, $y$ in $V(R)$
such that
\[
nd+b_1+x=nd+b_2+y\,.
\]
Hence, by
definition of the operation in $V(R)\sqcup W_{\sigma}^d(S_u)$, we get
that
\[
nd+b_1+\phi_u(x)=nd+b_2+\phi_u(y)\,.
\]
After restricting to the extreme
boundary of $S_u$ and using \cite[Lemma 3.4]{perjot} (and also that the
scale is finite), we see that
\[b_1+x=b_1+\phi_u(x)=b_2+\phi_u(y)=b_2+y\,.\]
Therefore $\ol{b_1}=\ol{b_2}$ in the quotient monoid, so we take
$c_1=c_2=0$ and $S_1=S_2=0$, and $x'=n\ol{d}$.

(ii). Let $u=[p]$, which is non-zero in $V(R)$, and continue to use the
notation of Theorem \ref{plasta} and Remark \ref{otia}. Assume that
$s,t\in d^{-1}(\{\infty\})\cap\partial_e S_u$ and that they are
different. As in \cite[Corollary 4.11]{percan}, there exist elements
$b_i$ in $W_{\sigma}^d(S_u)$ (for $i=1$, $2$) such that $b_1(s)=1$ and
${b_1}_{|\{s\}'}=d_{|\{s\}'}$, and that $b_2(t)=1$ and
${b_2}_{|\{t\}'}=d_{|\{t\}'}$. Observe that $b_1+d=2d=b_2+d$.

Assume, to reach a contradiction, that $\mulr/R$ is a $QB-$ring. Then
Remark \ref{otia} provides us with two order-ideals $S_i$ ($i=1,2$) of
$V(\mulr)$ and elements $\ol{c_i}$ in $\ol{S_i}$ such that $\ol{S_1}\cap
\ol{S_2}=0$ and $\ol{b_1}+\ol{c_1}=\ol{b_2}+\ol{c_2}$. If both $\ol{S_1}$
and $\ol{S_2}$ are non-zero, then $S_1$ and $S_2$ properly contain $V(R)$
so they must contain $V(L(R))$. But then
$\ol{S_1}\cap\ol{S_2}\supset\ol{V(L(R))}\ne 0$. We may assume without loss
of generality that $\ol{S_1}=0$, hence we have the equation
$\ol{b_1}=\ol{b_2}+\ol{c_2}$. We then find elements $x$ and $y$ in $V(R)$
such that $b_1+x=b_2+c_2+y$, that is, $b_1+\phi_u(x)=b_2+c_2+\phi_u(y)$.
Evaluating at $s$ we get the contradiction.
\end{proof}
\begin{prop}
Let $R$ be a ring in the class $\mathcal{N}$. Then $\mulr$ is a $QB-$ring
if and only if $R$ is elementary.
\end{prop}
\begin{proof}
Let $u$ be a non-zero element in $V(R)$, and set $S_u=\mathrm{St}(V(R),u)$
and $d=\sup\phi_u([p_n])$, where $(p_n)$ is a countable unit for $R$.
Assume that $R$ is non-elementary. Then Theorem \ref{plasta} ensures that
$V(\mathcal{M}(R))\cong V(R)\sqcup W_{\sigma}^d(S_u)$. Assume first that
$\mathcal{M}(R)$ is a $QB-$ring and that $R$ has finite scale. Then,
using the arguments in \cite[Section 7]{percan}, it follows that
$\mathcal{M}(R)$ is stably finite. Taking into account that
$\mathcal{M}(R)$ is moreover a prime ring, we get that
$\mathcal{M}(R)^{-1}_q=\mathcal{M}(R)^{-1}$. By \cite[Proposition
3.9]{qb}, $\mathcal{M}(R)$ has stable rank one, but this contradicts
Theorem \ref{plasta2}.

The scale is therefore not finite and so there exists at least one element
$s$ in $\partial_e S_u$ such that $d(s)=\infty$. Now condition (ii) in
Theorem \ref{mult1} implies that there can be only one infinite extremal
state, since otherwise $\mulr/R$ would not be a $QB-$ring, so neither
would be $\mulr$, in contradiction with our hypothesis.

Thus we have to deal with the case where \[V(\mathcal{M}(R))\cong
V(R)\sqcup\mathbb{R}^{++}\sqcup\{\infty\}\,.\] Let $\alpha>0$, take any
$a$ in $V(R)$ such that $\phi_u(a)>\alpha$, and consider the equation
$\infty+a=\infty=\infty+\alpha$. Another application of Remark \ref{otia}
provides us with elements $c_1$ and $c_2$ (one of which is zero, but not
both) such that $a+c_1=\alpha+c_2$. Now, since $a\in V(R)$, necessarily
$c_1\neq 0$ and so $c_2=0$. This implies that
$c_1+\phi_u(a)=\alpha>\phi_u(a)$, a contradiction.

Assume now that $R$ is elementary. We have already observed that in this
case $\mulr\cong\mathbb{B}(D)$ for a division ring $D$, and this is a
$QB-$ring.
\end{proof}
\begin{theor}
Let $R$ be a non-elementary ring in the class $\mathcal{N}$, and let $p$
be a non-zero idempotent in $R$. If $R$ has exactly one infinite extremal
state $s$ in $S_{[p]}$, then $\mulr/R$ is a $QB-$ring if and only if the
complementary face $\{s\}'$ is closed in $S_{[p]}$.
\end{theor}
\begin{proof}
Again, we let $u=[p]$ in $V(R)$ and adopt the notation of Theorem
\ref{plasta} and Remark \ref{otia}. Assume that $\{s\}'$ is closed and
that $nd+b_1+x=nd+b_2+y$, for some $x$, $y$ in $V(R)$ and some $n\geq 1$.
We may clearly assume that $b_i\notin V(R)$ for any $i$. We of course get
that $b_1+\phi_u(x)=b_2+\phi_u(y)$ in $\partial_e S_u\setminus\{s\}$. We
have to distinguish three cases:

If $b_1(s)=b_2(s)=\infty$, then $b_1+\phi_u(x)=b_2+\phi_u(y)$ in
$\partial_e S_u$, hence $b_1+x=b_1+y$ and therefore $\ol{b_1}=\ol{b_2}$
in $V(\mulr)/V(R)$. We let in this case $S_1=S_2=0$ and $x'=n\ol{d}$.

If $b_i(s)<\infty$ for all $i$, then without loss of generality we may
assume that the element $a=b_2(s)+\phi_u(y)(s)-(b_1(s)+\phi_u(x)(s))>0$.
Let $k$ be a natural number such that $k>a$. By the arguments in
\cite[Proposition 4.10]{percan}, there exists a function $e$ in
$\mathrm{LAff}(S_u)^{++}$ such that $e(s)=k-a$ and $e_{|\{s\}'}=k$. Now,
it is easy to check that
\[
b_1+\phi_u(x)+k=b_2+\phi_u(y)+e\,,
\]
whence
\[
b_1+ku+x=b_2+y+e\,.
\]
In the quotient monoid we therefore have $\ol{b_1}=\ol{b_2}+\ol{e}$.
Moreover, $e+d=ku+d$, and thus $\ol{e}+n\ol{d}=n\ol{d}$. We let in this
case $c_1=0$, $c_2=\ol{e}$, $x'=n\ol{d}$, $S_1=0$ and $S_2=V(\mulr)/V(R)$.

Finally, assume that $b_2(s)=\infty$ and $b_1(s)<\infty$. Let $k$ be a
natural number. Since $\{s\}'$ is closed, we can construct a function $e$
in $\mathrm{LAff}(S_u)^{++}$ such that $e(s)=\infty$ and $e_{|\{s\}'}=k$.
Now, $b_1+\phi_u(x)+e=b_2+\phi_u(y)+k$, and then $[b_1]+[e]=[b_2]$. As
before, we also have $e+d=ku+d$, so again $\ol{e}+n\ol{d}=n\ol{d}$. In
this case, $c_1=\ol{e}$, $c_2=0$, $x'=n\ol{d}$, $S_2=0$ and
$S_1=V(\mulr)/V(R)$.

To prove the converse direction, assume that $\mulr/R$ is a $QB-$ring and
that $\{s\}'$ is not closed. Then $F:=\ol{\{s\}'}\ne \{s\}'$, so we can
take an element $t$ in $F\setminus\{s\}'$. By \cite[Lemma 3.7]{perjot}
$F$ equals the closure of the convex hull $\mathrm{conv}(\partial_e
S_u\setminus\{s\})$. Therefore we can write $t=\lim s_n$, where
$s_n\in\mathrm{conv}(\partial_e S_u\setminus\{s\})$. On the other hand,
$S_u$ is the direct convex sum of $\{s\}$ and $\{s\}'$ (see \cite[Theorem
11.28]{poag}), so that there exist a real number $\alpha$ in $(0,1]$ and
an element $t_1$ in $\{s\}'$ such that
\[
t=\alpha s+(1-\alpha)t_1\,.
\]
By
using \cite[Corollary 4.11]{percan}, construct a function $d'$ such that
$d'(s)=1$ and $d'_{|\{s\}'}=d_{|\{s\}'}$, so $d+d'=d+d$. Since $\mulr/R$
is a $QB-$ring, there exist order-ideals $S_i$ ($i=1$, $2$) such that
$\ol{S_1}\cap \ol{S_2}=0$ in $V(\mulr)/V(R)$ and elements $\ol{c_i}$ in
$\ol{S_i}$ such that $\ol{d}+\ol{c_2}=\ol{d'}+\ol{c_1}$. Upstairs, this
means that
\[
d+c_2+x=d'+c_1+y
\]
for some $x$ and $y$ in $V(R)$. As in the proof of Theorem \ref{mult1},
we see that $\ol{S_1}=0$ or $\ol{S_2}=0$. Now $c_1(s)=\infty$ because
$d'(s)=1$, and we therefore conclude that $\ol{S_2}=0$, whence
$\ol{c_2}=0$ (and we may actually assume that $c_2=0$). The equation
\[
d+x=d'+c_1+y
\]
is equivalent to
\[
d+\phi_u(x)=d'+c_1+\phi_u(y)
\]
and restricting to the extreme boundary of
$\partial_e S_u$ we see that
\[
\phi_u(x)_{|\partial_e
S_u\setminus\{s\}}=(c_1+\phi_u(y))_{|\partial_e S_u\setminus\{s\}}\,.
\]
This implies (using the lower semi-continuity of $\phi_u(y)+c_1$) that
\[
\phi_u(x)(t)=\lim\limits_n
\phi_u(x)(s_n)=\lim\limits_n(c_1+\phi_u(y))(s_n)\geq
(c_1+\phi_u(y))(t)=\infty\,,
\]
which is plainly impossible.
\end{proof}

\end{document}